\documentclass[3p,preprint]{elsarticle_preprint}

\usepackage{amsmath}
\usepackage{amsfonts}
\usepackage{amssymb,latexsym}
\usepackage{amsthm}
\usepackage{graphicx}

\usepackage{color}
\usepackage{newlfont}
\usepackage{subfigure}
\usepackage{verbatim}
\usepackage{rotating}
\usepackage{multirow}
\usepackage{units}
\usepackage{bm}
\usepackage[english]{babel}
\usepackage[utf8]{inputenc}
\usepackage{algorithm}
\usepackage[noend]{algpseudocode}
\usepackage{booktabs}
\usepackage{caption}
\usepackage{hyperref}

\newcommand{\RR}{\mathbb{R}}

\newcommand{\mK}{\mathsf{K}}

\journal{R. Cavoretto}

\begin{document}

\begin{frontmatter}

\title{Adaptive LOOCV-based kernel methods for solving time-dependent BVPs}


\author[address-TO,address-GNCS]{Roberto Cavoretto} 
\ead{roberto.cavoretto@unito.it}



\address[address-TO]{Department of Mathematics \lq\lq Giuseppe Peano\rq\rq, University of Torino, via Carlo Alberto 10, 10123 Torino, Italy}

\address[address-GNCS]{Member of the INdAM Research group GNCS}

\begin{abstract}

In this paper we propose an adaptive scheme for the solution of time-dependent boundary value problems (BVPs). To solve numerically these problems, we consider the kernel-based method of lines that allows us to split the spatial and time derivatives, dealing with each separately. This adaptive algorithm is based on a leave-one-out cross validation (LOOCV) technique, which is employed as an error indicator. By this scheme, we can first detect the domain areas where the error is estimated to be too large -- generally due to steep variations or quick changes in the solution -- and then accordingly enhance the numerical solution by applying a two-point refinement strategy. Numerical experiments show the efficacy and performance of our adaptive refinement method.  

\end{abstract}

\begin{keyword}
meshless approximation\sep kernel methods\sep adaptive algorithms\sep refinement schemes\sep partial differential equations
\MSC[2020] 65D12, 65D15, 65M20
\end{keyword}


\end{frontmatter}


\section{Introduction} \label{sec:1}


In the works \cite{cav20a,cav20b} adaptive refinement schemes for solving elliptic boundary value problems (BVPs) by radial basis function collocation methods are proposed. In particular, the work \cite{cav20b} provides an efficient algorithm to adaptively compute the numerical solution of differential problems defined over 2D irregular domains. It is based on a leave-one-out cross validation (LOOCV) error indicator, which is combined with a specific refinement strategy. This adaptive method can successfully be applied to time-independent problems, while the solution of time-dependent BVPs is still an open problem. Note that the above-mentioned LOOCV refers to Rippa's formula \cite{rip99}, which was often used in applied mathematics to optimally compute the kernel shape parameter (see, e.g., \cite{cav21b,cav22,lin21,mar21,yao15}, and \cite[Chapter 14]{fas15} for an overview).

In this article we thus extend the work \cite{cav20b} to solve time-dependent BVPs. By doing so, we choose to use the kernel-based method of lines, which enables us to distinguish between spatial and physical components, dealing with each separately (see \cite{che07,fas15}). The adaptive process is characterized by an iterated use of the LOOCV-based error indicator. The latter is applied as a predictor to estimate at any time level whether the approximate solution found by the kernel-based method of lines is accurate enough, or instead the error is over a prescribed threshold. Specifically, when the refinement process needs to be employed, we apply iteratively a two-point refinement procedure.  As an example, we test the resulting adaptive scheme on some nonlinear Burgers and Allen-Cahn BVPs (see \cite{dri07,hua94,yan19}), whose solutions may sometimes exhibit steep gradients, corners and rapid changes. These solutions, if exist, are generally difficult to obtain analytically, so a numerical approach must necessarily be used \cite{naq17}. Efficacy and performance of our adaptive refinement procedure are shown by our numerical experiments. 


The paper is organized as follows. In Section \ref{sec:2} we introduce the kernel-based method of lines for BVPs. Section \ref{sec:3} presents the adaptive refinement algorithm together with its LOOCV-based error indicator. In Section \ref{sec:4} we show some numerical results so as to illustrate the performance of our adaptive scheme. 

\section{Kernel-based method of lines for BVPs} \label{sec:2}

Given a domain $\Omega = [a,b] \subset \RR$, we consider a time-dependent BVP problem of the form
\begin{eqnarray} 
\begin{array}{rl}
\label{PDE-BC}
\displaystyle{\frac{\partial}{\partial t}}u(x,t)={\cal L}u(x,t)+f(x,t),& \qquad x\in \Omega, \,\, t>t_0, \medskip \\
0={\cal B}u(x,t)+g(x,t), & \qquad x\in \partial\Omega, \,\, t \geq t_0, \medskip \\
u(x,t)=u_0(x), & \qquad x\in \Omega, \,\, t = t_0,
\end{array}
\end{eqnarray}
where ${\cal L}$ is a linear or nonlinear partial differential operator, and ${\cal B}$ is a boundary operator. The functions $f$, $g$ and $u_0$ are provided as data to the problem. We assume that the kernel centers $Z_N=\{z_1,\ldots,z_N\} \subseteq \Omega$ and collocation points $X_N=\{x_1,\ldots,x_N\}$ coincide, i.e. $Z_N \equiv X_N$. Then, supposing the collocation points are ordered such that $x_1 < x_2 < \cdots < x_{N-1} < x_{N}$, the set $X_N$ is subdivided into the subsets $X_{N_I}=\{x_2,\ldots,x_{N-1}\}$ of interior points and $X_{N_B}=\{x_{1},x_N\}$ of boundary points, so that $X_N = X_{N_I} \cup X_{N_B}$, where the subscript $N_I$ and $N_B$ denote the number of points in the interior and boundary, respectively. In this case, we fix $a \equiv x_{1}$ and $b \equiv x_N$.


For the method of lines we express the approximate solution $\hat{u}$ as a linear combination of radial kernel $K:\Omega\times \Omega \rightarrow \RR$, i.e.,
\begin{eqnarray} \label{hatu}
\hat{u}(x,t) = \sum_{j=1}^{N} c_j(t) K(x,x_j) = \boldsymbol{k}(x)^T\boldsymbol{c}(t),
\end{eqnarray}
where $\boldsymbol{k}(x)^T=(K(x,x_1) \cdots K(x,x_N))$ and $\boldsymbol{c}(t)=(\boldsymbol{c}_1(t) \cdots \boldsymbol{c}_N(t))^T$ are unknown coefficients to be found. Associated with the kernel $K$ we can define a function of a single scalar variable $\kappa:\RR_0^+\rightarrow \RR$ depending on a \textsl{shape parameter} $\varepsilon > 0$ such that $K(x,x_j)=\kappa_{\varepsilon}(r)=\kappa(\varepsilon r)$, with $r=|x-x_j|$ denoting the Euclidean distance between $x$ and $x_j$ in $\RR$. Though the choice of $\varepsilon$ often does not depend on the point $x_j$, for our adaptive purposes we will consider a variable selection $\varepsilon_j$ based on the spacing between $x$ and $x_j$. Accordingly, the kernel will be more localizing or picked when the distance between two points is small, while it will turn out to be less localizing or flatter in case of a large distance between the points. More precisely, starting from an initial shape parameter $\varepsilon_0$, the variable shape parameter $\varepsilon_j$ is defined as:
\begin{eqnarray} 	\label{epsds}
\varepsilon_j = \left\{
\begin{array}{ll}
\displaystyle{\varepsilon_0\frac{1}{|x_{j} - x_{j+1}|}}, & \qquad j = 1 , \medskip \\
\displaystyle{\varepsilon_0\min \left\{\frac{1}{|x_{j-1} - x_{j}|}, \frac{1}{|x_j - x_{j+1}|}\right\}}, & \qquad j=2,\ldots,N-1, \medskip \\
\displaystyle{\varepsilon_0\frac{1}{|x_{j-1} - x_{j}|}}, & \qquad j = N.
\end{array}
\right.
\end{eqnarray}

The first step to get a method of lines solution is to compute a kernel-based interpolant, by evaluating our kernel solution \eqref{hatu} at fixed time $t=0$. Thus, requiring to solve the linear system
    \begin{align} \label{sysIC}
\left(
\begin{array}{c}
\boldsymbol{k}(x_1)^T \\
\boldsymbol{k}(x_2)^T \\
\vdots \\
\boldsymbol{k}(x_{N-1})^T \\
\boldsymbol{k}(x_{N})^T       
\end{array}
\right)
\boldsymbol{c}(0)
=
\left(
\begin{array}{c}
u_0(x_1)\\
u_0(x_2)\\
\vdots\\
u_0(x_{N-1})\\
u_0(x_{N})
\end{array}
\right),
  \end{align}
we obtain the $N$ initial conditions $\boldsymbol{c}(0)$. Then, to create a system of $N$ ordinary differential equations for the coefficient expantion $\boldsymbol{c}(t)$, we can semi-discretize the PDE in \eqref{PDE-BC} by collocating at the $N-2$ interior centers and generating the $(N-2) \times N$ system
    \begin{align} \label{sysPDE}
\left(
\begin{array}{c}
\boldsymbol{k}(x_2)^T \\
\vdots \\
\boldsymbol{k}(x_{N-1})^T
\end{array}
\right)
\displaystyle{\frac{\partial}{\partial t}} \boldsymbol{c}(t)
=
\left(
\begin{array}{c}
{\cal L}\boldsymbol{k}(x_2)^T \\
\vdots \\
{\cal L}\boldsymbol{k}(x_{N-1})^T
\end{array}
\right)
\boldsymbol{c}(t)
+
\left(
\begin{array}{c}
f(x_2,t)\\
\vdots\\
f(x_{N-1},t)
\end{array}
\right).
  \end{align}
Discretizing the boundary conditions in \eqref{PDE-BC}, we derive the $2 \times N$ system
\begin{align} \label{sysBC}
\left(
\begin{array}{c}
0 \\
0
\end{array}
\right)
=
\left(
\begin{array}{c}
{\cal B}\boldsymbol{k}(x_{1})^T \\
{\cal B}\boldsymbol{k}(x_{N})^T
\end{array}
\right)
\boldsymbol{c}(t)
+
\left(
\begin{array}{c}
g(x_{1},t)\\
g(x_{N},t)
\end{array}
\right).
  \end{align}
Therefore, coupling \eqref{sysPDE} and \eqref{sysBC}, we yield a full set of $N$ differential algebraic equations (DAEs) of the form
    \begin{align} \label{sysPDE+BC}
\left(
\begin{array}{c}
0 \\
\boldsymbol{k}(x_1)^T \\
\vdots \\
\boldsymbol{k}(x_{N-1})^T\\
0
\end{array}
\right)
\displaystyle{\frac{\partial}{\partial t}} \boldsymbol{c}(t)
=
\left(
\begin{array}{c}
{\cal B}\boldsymbol{k}(x_1)^T \\
{\cal L}\boldsymbol{k}(x_2)^T \\
\vdots \\
{\cal L}\boldsymbol{k}(x_{N-1})^T\\
{\cal B}\boldsymbol{k}(x_{N})^T
\end{array}
\right)
\boldsymbol{c}(t)
+
\left(
\begin{array}{c}
g(x_{1},t)\\
f(x_2,t)\\
\vdots\\
f(x_{N-1},t)\\
g(x_{N},t)
\end{array}
\right).
  \end{align}
As a result, the system of DAEs \eqref{sysPDE+BC} together with the initial conditions derived from \eqref{sysIC} allows us to find the kernel-based method of lines solution. Thus, if we discretize in time, at each time step $t=t_n$, with $n=1,\ldots,M$, where $M$ denote the maximum level of time integration, the solution vector $\boldsymbol{c}(t_n)=(\boldsymbol{c}_1(t_n) \cdots \boldsymbol{c}_N(t_n))^T$ is expressed as $\boldsymbol{c}_n=(\boldsymbol{c}_{n1} \cdots \boldsymbol{c}_{nN})^T$. Similarly, we denote the approximate solution vector of \eqref{hatu} as $\hat{\boldsymbol{u}}_n=(\hat{u}_{n1} \cdots \hat{u}_{nN})^T$, i.e., $\hat{u}_{nj}=\hat{u}(x_j,t_n)$, for any $j=1,\ldots,N$ and $n=1,\ldots,M$.

\section{LOOCV-based adaptive algorithm} \label{sec:3}

\subsection{Error indicator}

Firsty, we introduce a kernel-based interpolant on the given data $U_N=\{(x_j,\hat{u}_n(x_j))\}_{j=1}^N$ of the form      
\begin{eqnarray} \label{interps}
s_n(x) = \sum_{j=1}^{N} \alpha_{nj} K(x,x_j) = \boldsymbol{k}(x)^T\boldsymbol{\alpha}_n,
\end{eqnarray}
where $\boldsymbol{k}(x)^T$ is defined as in \eqref{hatu}, the subscript $n$ denotes a generic time level, and the coefficients $\boldsymbol{\alpha}_n=(\alpha_{n1} \cdots \alpha_{nN})^T$ are obtained by solving the linear system
\begin{align} \label{sysK}
\mK \boldsymbol{\alpha}_n = \hat{\boldsymbol{u}}_n,
\end{align}
where $\mK_{ij} = (K(x_i,x_j))_{i,j=1}^N$ and $\hat{\boldsymbol{u}}_n=(u_{n1} \cdots u_{nN})^T$.
For the sake of simplicity, in \eqref{interps} we omit to employ the symbol $n$ on the spatial data $x$ or $x_j$, although in an adaptive scheme they are obviously time-dependent.

In order to make an estimate of the error on the approximate solution and detect which domain areas are to be refined by the addition of new discretization points, we split the data set $U_N$ into two parts: i) a \emph{training} data set $T_N=\{(x_j,\hat{u}_n(x_j))\}_{j=1,j\neq k}^N$ consisting of $N - 1$ data to obtain a \lq\lq partial interpolation\rq\rq; ii) a \emph{validation} data set $V_N=\{(x_j,\hat{u}_n(x_j))\}_{j=k}$ that contains a single (remaining) data used to compute the error. Repeating in turn this procedure by varying the index $k$, i.e., for each of the $N$ given data, we obtain a vector of error estimates that enables us to identify an error indicator based on a LOOCV strategy \cite{fas07b}.

To simplify the description of LOOCV method, we define by
\begin{equation} \label{defvec}
\boldsymbol{x}^{[k]}=(x_1 \cdots x_{k-1} x_{k+1} \cdots x_N)^T, \qquad \boldsymbol{\hat{u}}_n^{[k]}=(\hat{u}_{n1} \cdots,\hat{u}_{nk-1} \hat{u}_{nk+1} \cdots \hat{u}_{nN})^T,
\end{equation}
the vectors of interpolation points and corresponding values derived from method of lines solution, where the superscript $[k]$ denotes the lack of removed data $(x_{k},\hat{u}_{nk})$. In the following, all other qualities are represented similarly.  

Now, for a fixed index $k \in \{1,\ldots,N\}$ and a fixed value of the shape parameter $\varepsilon$, we define the partial kernel-based interpolant
\begin{align} \label{pinterp}
	s_n^{[k]}(x) = \sum_{j=1,\ j \neq k}^{N} \alpha_{nj} K (x,x_j),
\end{align}
whose coefficients $\alpha_{nj}$ are determined by interpolating the training data set, i.e. $s_n^{[k]}(x_i) = \hat{u}_n(x_i)$, $i=1,\ldots,k-1,k+1,\ldots,N$. This is equivalent to solving the $(N - 1) \times (N - 1)$ linear system
\begin{align} \label{sysinterp}
\mK^{[k]} \boldsymbol{\alpha}_{n}^{[k]} = \boldsymbol{\hat{u}}_n^{[k]},
\end{align}
where $\mK^{[k]}$ is obtained from the full interpolation matrix $\mK$ in \eqref{sysK}, with entries $(\mK)_{ij}=(K(x_i,x_j))_{i,j=1}^N$, by removing the $k$th row and the $k$th column, while $\boldsymbol{\alpha}_{n}^{[k]}=(\alpha_{n1} \cdots \alpha_{nk-1} \alpha_{nk+1} \cdots \alpha_{nN})^T$ and $\boldsymbol{\hat{u}}_n^{[k]}$ is given in the right hand side of \eqref{defvec}. From \eqref{sysinterp} we can note that the kernel is here used as a predictor to interpolate solution values at the collocation points. This stage is particularly relevant for the construction of an adaptive scheme, because it permits to identify the domain regions where adding \lq\lq new\rq\rq\ collocation points.     

Now, instead of solving the system \eqref{sysinterp}, as suggested in \cite{rip99}, we can make a prediction on the behavior of the numerical solution by computing the absolute error components $e_{nk}$ by the rule 
\begin{align} \label{ekeps}
e_{nk} = \left| \frac{\alpha_{nk}}{\mK_{kk}^{-1}} \right|, \qquad k=1,\ldots,N,
\end{align}
where $\alpha_{nk}$ is the $k$th coefficient of the full kernel-based interpolant \eqref{interps} and $\mK_{kk}^{-1}$ is the $k$th diagonal element of the $N \times N$ interpolation matrix $\mK^{-1}$. By applying the rule \eqref{ekeps}, the computational complexity is ${\cal O}(N^3)$. As a consequence, we obtain a vector $\boldsymbol{e}_{n} = (e_{n1},\ldots,e_{nN})^T$, which can be used as an error indicator to identify the domain areas that need a refinement in proximity of point $x_k$.



\subsection{Iterated refinement scheme} \label{sec:32}

In this subsection we describe the iterative scheme based on the LOOCV error indicator. Since we use the method of lines for solving time-dependent BVPs, we apply our adaptive refinement algorithm at fixed time $t_n$, for $n=0,1,\ldots,M$. In this way, once fixed the time component, the adaptive scheme works iteratively on the spatial term, for any $n$. Notice that, although in our iterative procedure the collocation points change at any time level $t_n$, we imply their dependence from $n$ but, for the sake of simplicity, we neglect it in this discussion.

At first, we thus define an initial set $X_{N^{(1)}}= \{ x_1^{(1)}, \ldots, x_{N^{(1)}}^{(1)}\}$ of equally-spaced collocation points, with $X_{N^{(1)}}\equiv X_N$. It is then split into two subsets, i.e. the set $X_{N_I^{(i)}}=\{ x_{2}^{(i)}, \ldots, x_{N^{(i)}-1}^{(i)}\}$ of interior points, and the set $X_{N_B^{(i)}}=\{ x_{1}^{(i)}, x_{N^{(i)}}^{(i)}\}$ of boundary points, with $X_{N_B^{(i)}} \equiv X_{N_B}$ for any $i=1,2,\ldots$, where the subscript $_{N^{(i)}}$ denotes the number of points in the $i$th iteration at any time level $n$. 

Then, fixed a tolerance $\tau >0$, from \eqref{ekeps} we can define an iterated LOOCV-based error indicator as follows
\begin{align} \label{ekepsi}
e_{nk}^{(i)} = \left| \frac{\alpha_{nk}}{\mK_{kk}^{-1}} \right|, \qquad k=1,\ldots,N^{(i)}, \qquad  i=1,2,\ldots,
\end{align}
where in $|\cdot|$ some reference to the iteration is overlooked to simplify the notation.

Now, if the error indicator \eqref{ekepsi} is such that $e_{nk}^{(i)} > \tau$, a refinement is applied nearby $x_k$. By doing so, we first compute the so-called \textsl{separation distance} 
\begin{eqnarray} 	\label{sd}
q_{X_{N^{(i)}}} = \frac{1}{2}\min_{u \neq v} |x_u^{(i)} - x_v^{(i)}|, \qquad x_u^{(i)}, x_v^{(i)} \in X_{N^{(i)}}, \qquad i=1,2,\ldots,
\end{eqnarray}
and then sum up or subtract the quantity in \eqref{sd} to the point $x_k$. Specifically, the refinement procedure permits the addition of two points close to $x_k$, i.e., $x_k+q_{X_{N^{(i)}}}$ and $x_k-q_{X_{N^{(i)}}}$. Here, in general, we have to remark that in the refinement process multiple collocation points could be created. Thus, in order to avoid this issue, our procedure makes a check and, if necessary, excludes duplicated points. The set $X_{N^{(i)}}$ is then updated. The choice of using a two-point approach is probably the simplest one, but it allows us a precise refinement in specific parts of the integration interval $\Omega$.

Finally, when all components of error terms in \eqref{ekepsi} are less than or equal to $\tau$, the refinement process stops and the final set of collocation points is returned.
\section{Numerical results} \label{sec:4}

In this section we present numerical results obtained by solving some time-dependent BVPs. By doing so, we validate our adaptive refinement scheme by considering a few examples of nonlinear Burgers and Allen-Cahn equations defined on a space-time domain $(a,b) \times (0,T]$. All the experiments are run in \textsc{Matlab} and run on a laptop with Intel(R) Core(TM) i7-6500U CPU 2.50 GHz processor and 8GB RAM.

In our tests we show the results obtained by applying the method of lines and using as a basis function in \eqref{hatu} the multiquadric kernel $\kappa(\varepsilon_j r) = \left(1+\varepsilon_j^2 r^2\right)^{1/2}$, where $\varepsilon_j$ is the shape parameter chosen as in \eqref{epsds}. Then, for the integration in time we use the \textsc{Matlab} \texttt{ode15s} built-in routine and we assess the accuracy of the adaptive method by computing the \emph{root mean square error} at time $t$ (shortly, RMSE($t$)), i.e.,
\begin{align*} 
	\mbox{RMSE}(t) = \left(\frac{1}{N_{fin}}\sum_{i=1}^{N_{fin}} |u(x_i,t) - \hat{u}(x_i,t)|^2\right)^{1/2},
\end{align*}
where $x_i$, $i=1,\ldots,N_{fin}$, are the final points at fixed time $t=t_n$, for $n=0,1,\ldots,M$. Moreover, since the quality of predictions due to our LOOCV-based error indicator is dependent from matrix $\mK$ in \eqref{sysK}, we compute its \emph{condition number} (CN($t$)) by means of the \textsc{Matlab} command \texttt{cond}. The choice of various parameters are obviously problem dependent, so they will be defined case-by-case in the following. Finally, as a measure of computational efficiency we report the CPU time (CPU($t$)) expressed in seconds.


\subsection{Case study 1: Solution of Burgers equations} \label{sec:41}

In this subsection we consider time-dependent BVPs of the form \eqref{PDE-BC}. More specifically, in this study we take the following nonlinear Burgers equation \cite{bur48} 
\begin{align} \label{burgers}
\frac{\partial}{\partial t} u(x,t) = \nu \frac{\partial^2}{\partial x^2}u(x,t) - u(x,t) \frac{\partial}{\partial x} u(x,t), \qquad x \in (a,b), \quad t \in (t_0,T],
\end{align}
where the given parameter $\nu > 0$ develops a shock, and the solution exhibits moving fronts that can be made arbitrarily sharp by decreasing the kinematic viscosity $\nu$. The Burgers equation finds applications in several fields, for instance including modelling of fluid dynamic, shock wave formation, traffic flow and turbulence (see \cite{naq17}).

\begin{figure}[ht!]
\centering
{
\includegraphics[scale=0.26]{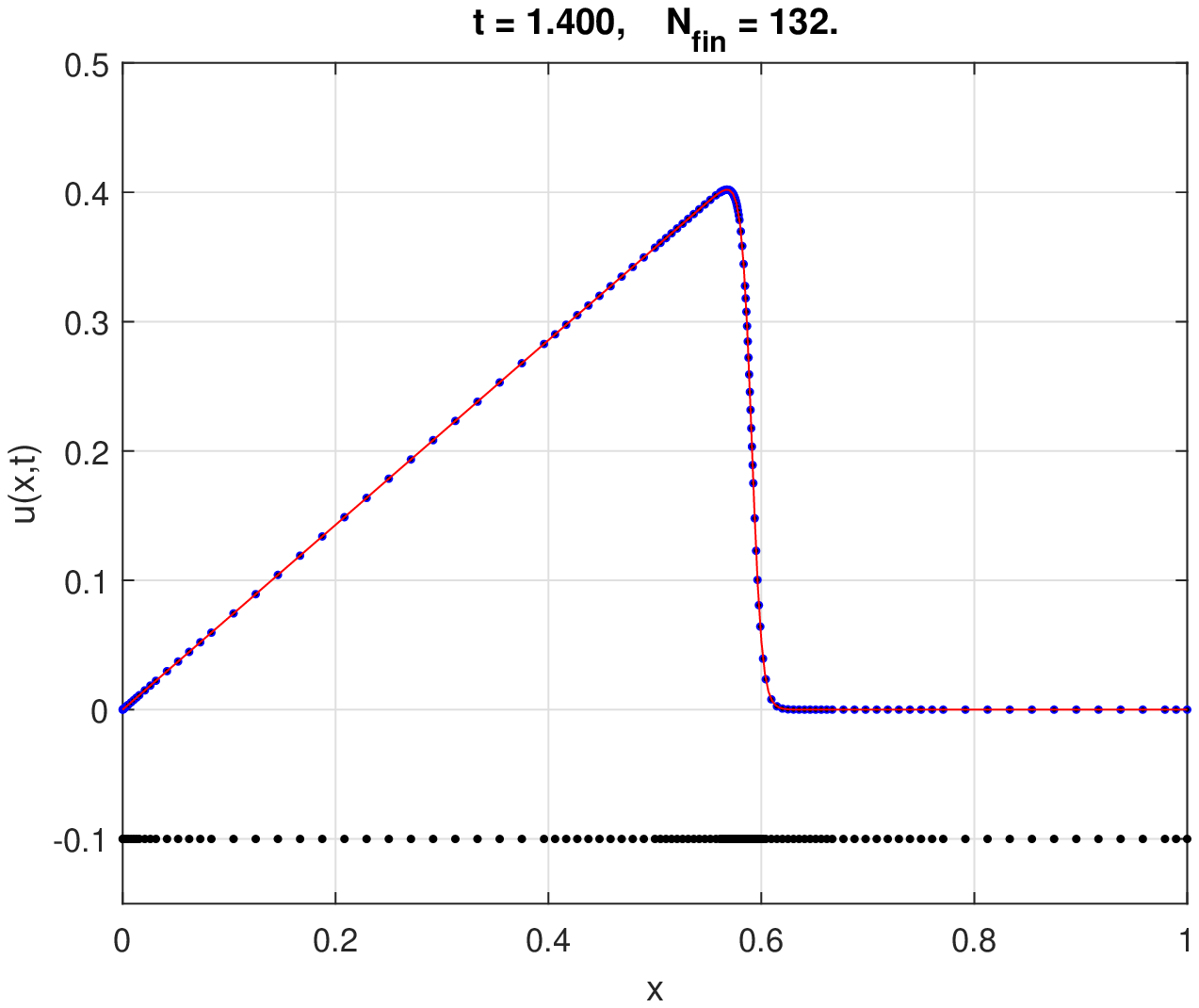} 
\includegraphics[scale=0.26]{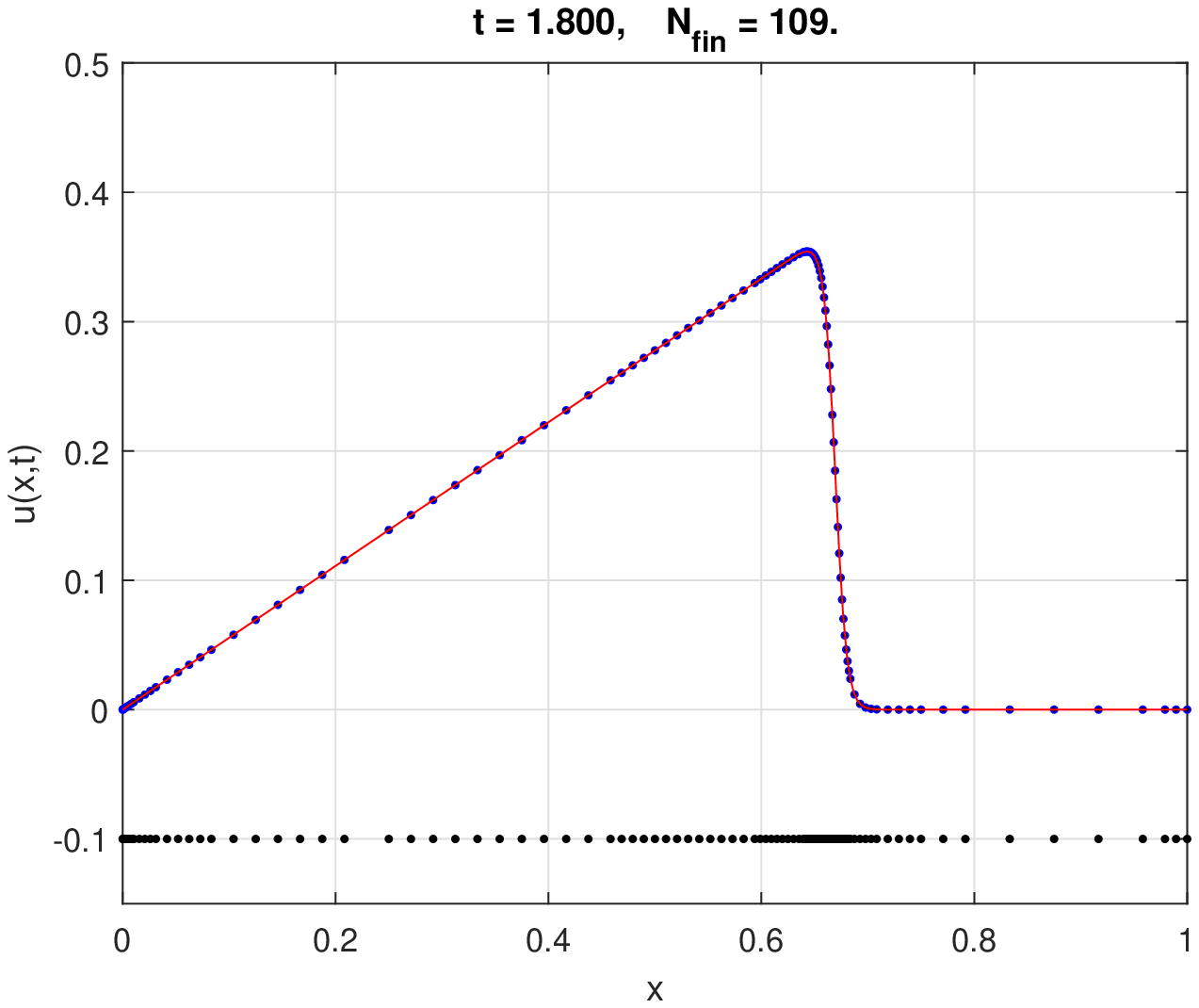} 
\includegraphics[scale=0.26]{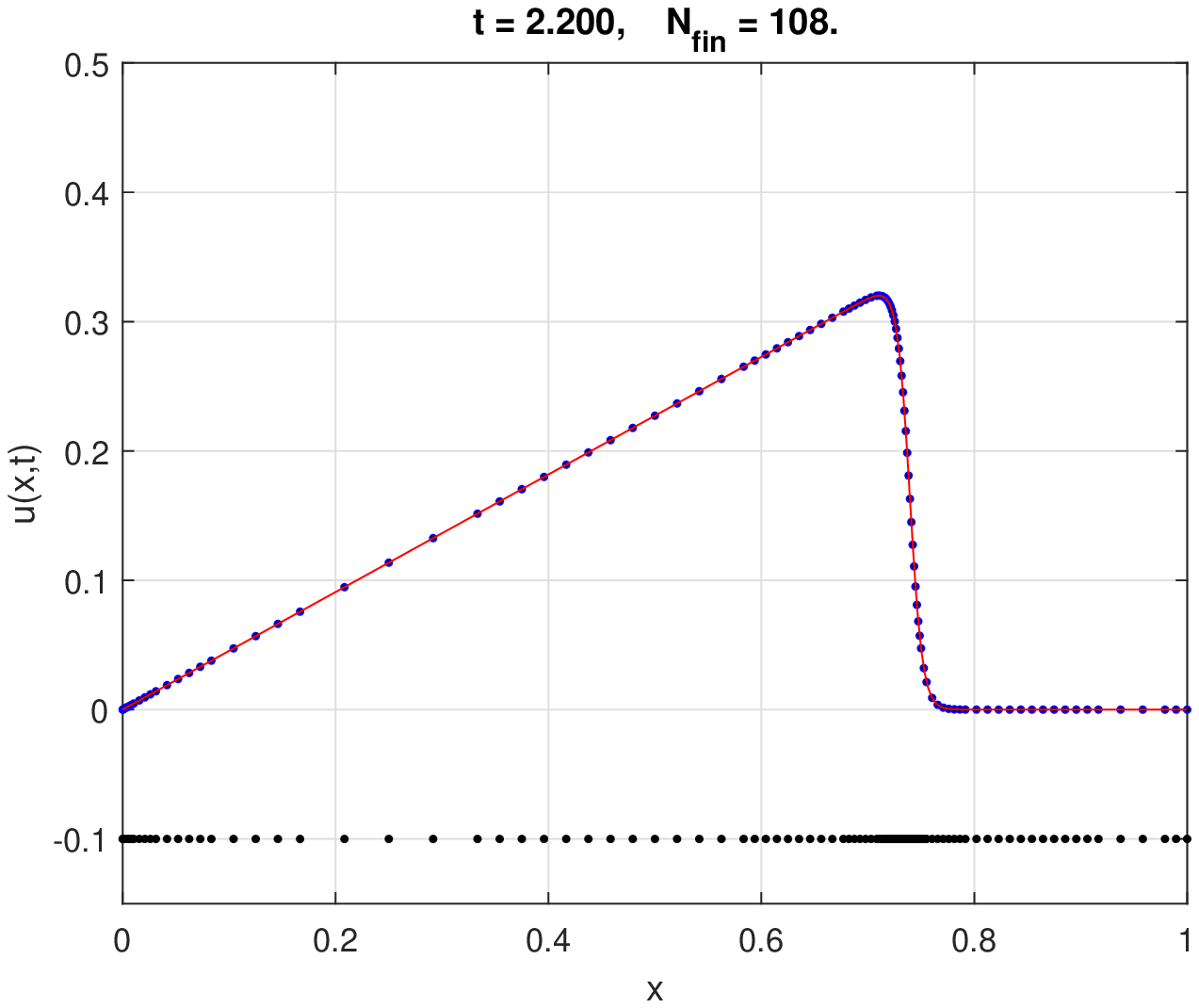} 
\includegraphics[scale=0.26]{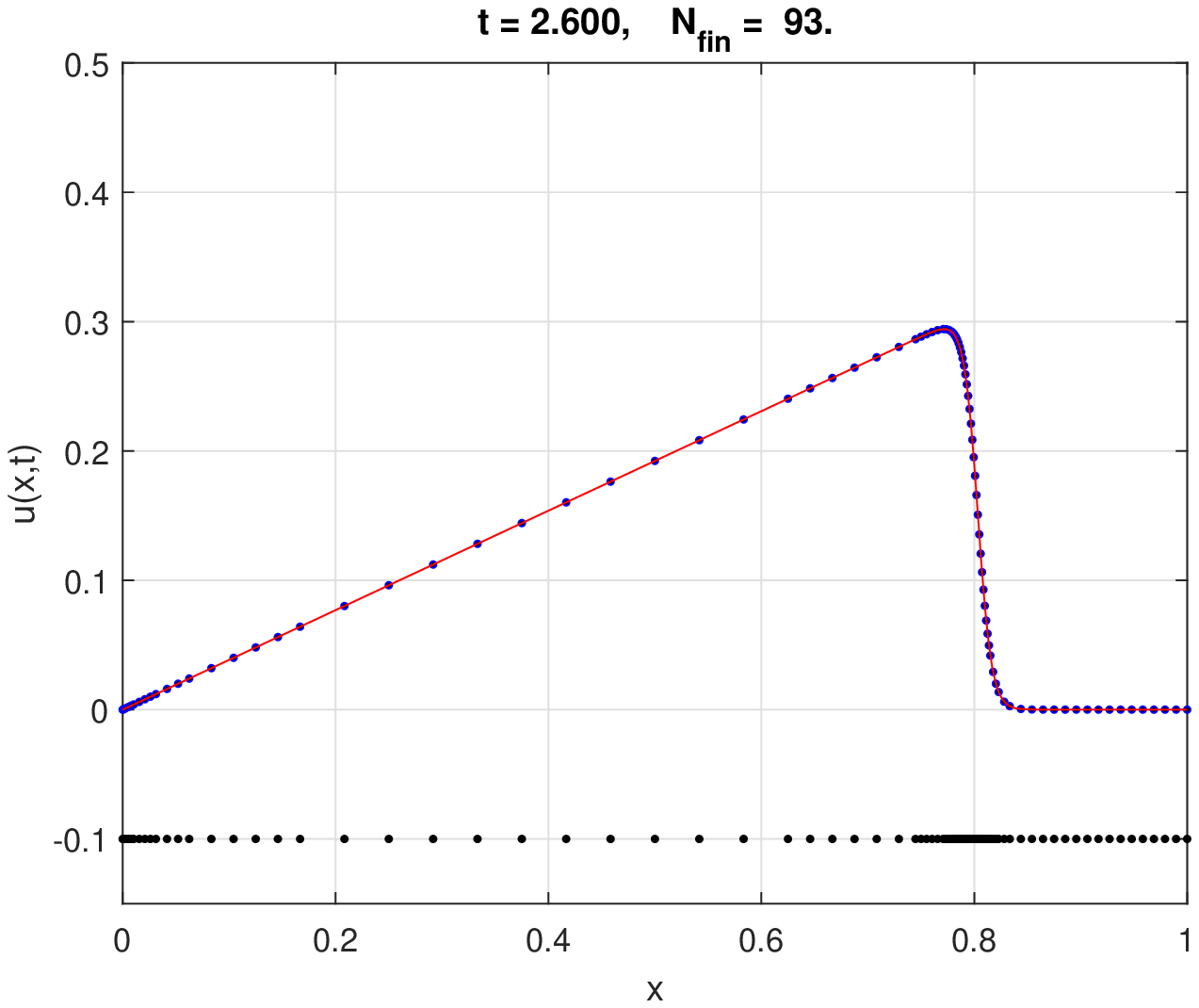}
}
\caption{Approximate solutions for Burgers equation of Example 1: $t=1.4$, $t=1.8$, $t=2.2$ and $t=2.6$ (left to right). Such graphs refer to results given in Table \ref{tab:1}.}
\label{fig:1}
\end{figure}

\vskip 0.2cm
\noindent \textbf{Example 1.} We consider the shock-like solution of \eqref{burgers}, whose analytical solution given in \cite{guo16,har96} is
\begin{align} \label{ex1_solburg}
\displaystyle{u(x,t) = \frac{\frac{x}{t}}{1+\sqrt{\frac{t}{c}}\exp\left(\frac{x^2}{4\nu t}\right)}}, \qquad t\in [t_0,T],
\end{align} 
where $c=\exp\left(\frac{1}{8\nu}\right)$. The initial condition for this example is obtained from the exact solution \eqref{ex1_solburg} when $t = t_0$ is used, while the boundary conditions are
\begin{align} \label{bc_burg}
u(a,t) = u(b,t) = 0, \qquad t\in [t_0,T].
\end{align}
By setting $a=1$, $b=1$, $t_0=1$, $T=3$, with the parameters $\nu = 10^{-3}$ and $\varepsilon_0=0.75$, the method starts initially with $N=13$ equally-spaced collocation points in the spatial interval $[0,1]$ and the adaptation occurs at every time level $t_n$, with $n=1,\ldots,M$. Here, we assume that the temporal interval $[1,3]$ is discretized by fixing $M=51$. In Table \ref{tab:1} we show the numerical results obtained for Burgers equation at different time levels, i.e., $t=1.4, 1.8, 2.2, 2.6, 3.0$. In such cases, we fix the prescribed threshold value $\tau=10^{-4}$ as a stop criterion of the LOOCV-based error indicator in \eqref{ekepsi}. From this table we can observe that the adaptive algorithm converges quite quickly in a final number $N_{fin}$ of points. Furthermore, the variable selection of the kernel shape parameters results in accurate results, also providing a control of the conditioning. In Figure \ref{fig:1} we depict the numerical solutions together with the final configurations of collocation points for some values of $t$. Note that the adaptive refinement scheme is able to detect and increase the points in the domain areas where steep variations of the solution are present.

\begin{table}[ht!]
{\small
\begin{center} 
{
\begin{tabular}{cccccc} 
\midrule
 $t$ & $N_{fin}$ & RMSE($t$) & CN($t$) & CPU($t$)   \\ \midrule
 $1.4$ & $132$  & $1.2$e${-4}$ & $6.8$e${+6}$ & $0.05$  \\
 $1.8$ & $109$  & $2.0$e${-4}$ & $6.9$e${+6}$ & $0.04$  \\
 $2.2$ & $108$  & $2.6$e${-4}$ & $7.3$e${+6}$ & $0.04$  \\
 $2.6$ & $93$  & $2.9$e${-4}$ & $3.7$e${+6}$ & $0.03$  \\
 $3.0$ & $98$  & $3.0$e${-4}$ & $3.9$e${+6}$ & $0.02$  \\
\midrule
\end{tabular}
}
\caption{Results obtained for Burgers equation of Example 1 by using the threshold value $\tau = 10^{-4}$ at some time levels $t$.} 
\label{tab:1}
\end{center}
}
\end{table}
 
While in Table \ref{tab:1} the threshold value $\tau$ is fixed and $t$ varies, in Table \ref{tab:2} we report the results for different values of $\tau$ going from $10^{-2}$ to $10^{-5}$, at fixed time $t=T=3$. From these experiments it is evident that the precision of the numerical scheme -- as well as the number of final points required by the adaptive process -- increases for smaller and more demanding values of $\tau$. However, even when the threshold is quite demanding, the refinement algorithm turns out to be efficient in terms of CPU times. 

\begin{table}[ht!]
{\small
\begin{center} 
{
\begin{tabular}{cccccc} 
\midrule
 $\tau$ & $N_{fin}$ & RMSE($t$) & CN($t$) & CPU($t$)   \\ \midrule
 $10^{-2}$ & $34$  & $6.0$e${-4}$ & $1.4$e${+5}$ & $<0.01$  \\
 $10^{-3}$ & $53$  & $3.6$e${-4}$ & $4.6$e${+5}$ & $0.01$  \\
 $10^{-4}$ & $98$  & $3.0$e${-4}$ & $3.9$e${+6}$ & $0.02$  \\
 $10^{-5}$ & $193$  & $7.4$e${-5}$ & $1.5$e${+8}$ & $0.08$  \\
\midrule
\end{tabular}
}
\caption{Results obtained for Burgers equation of Example 1 by using various values of threshold $\tau$ at time $t=T=3$.} 
\label{tab:2}
\end{center}
}
\end{table}


\vskip 0.2cm
\noindent \textbf{Example 2.} We consider the moving front problem given by the Burgers equation \eqref{burgers}, which is known as a common test problem for adaptive methods (see \cite{dri07,hua94}). The initial condition at time $t_0=0$ is
\begin{align*} 
u(x,0)=\sin(2\pi x) + \frac{1}{2} \sin(\pi x),
\end{align*}
while the boundary conditions are defined as in \eqref{bc_burg}, with $a=0$, $b=1$ and $T=1$. Taking homogeneous boundary conditions leads to a reduction of the wave amplitude when increasing time. The solution thus represents a wave that creates a steep front of width ${\cal O}(\nu)$ moving towards $x=1$. This is one of the main reasons for which this test BVP is quite stringent. Indeed, the solution steepens with increasing time $t$, thus becoming difficult to solve spatially \cite{naq17}.

\begin{figure}[ht!]
\centering
{
\includegraphics[scale=0.26]{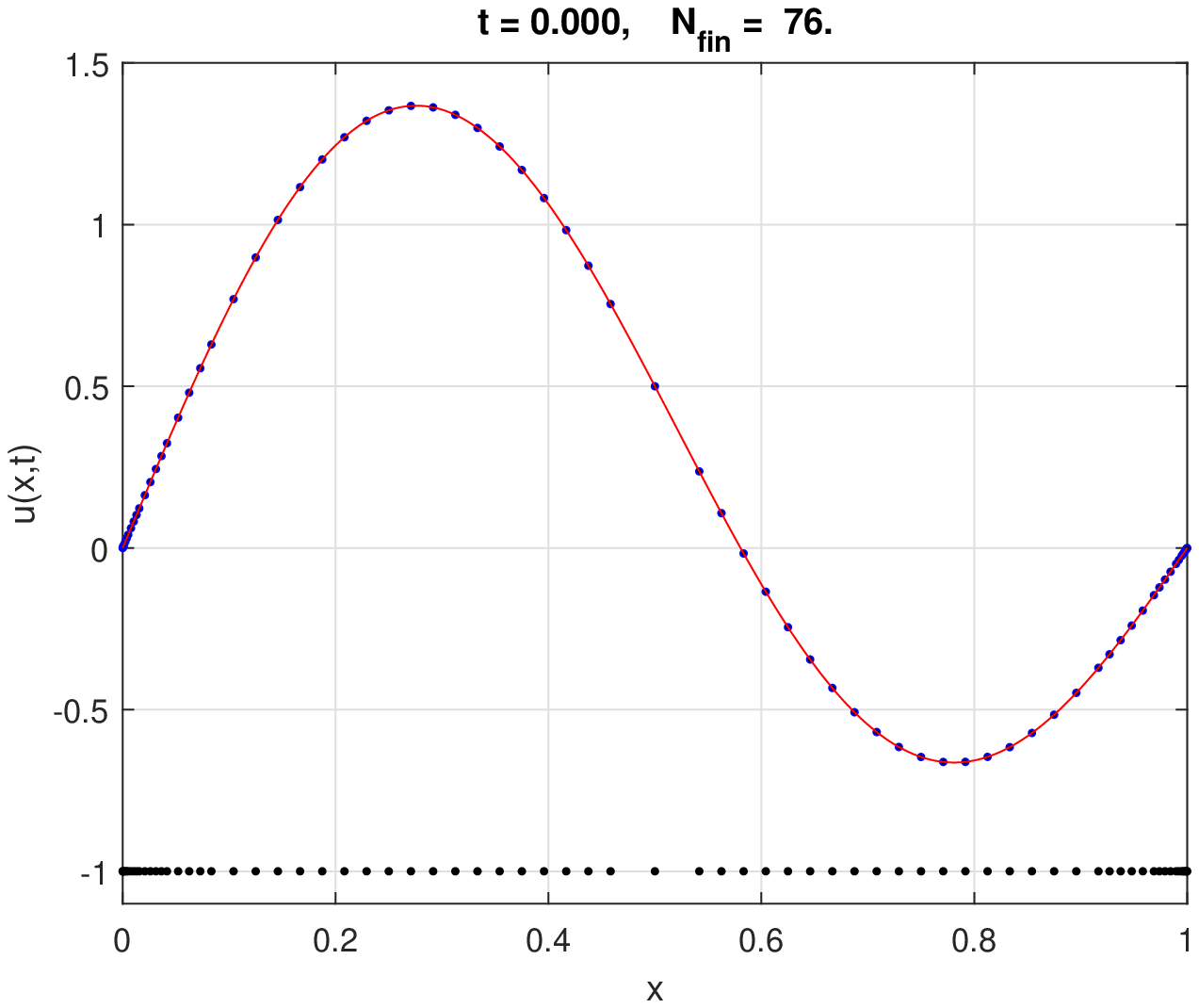} 
\includegraphics[scale=0.26]{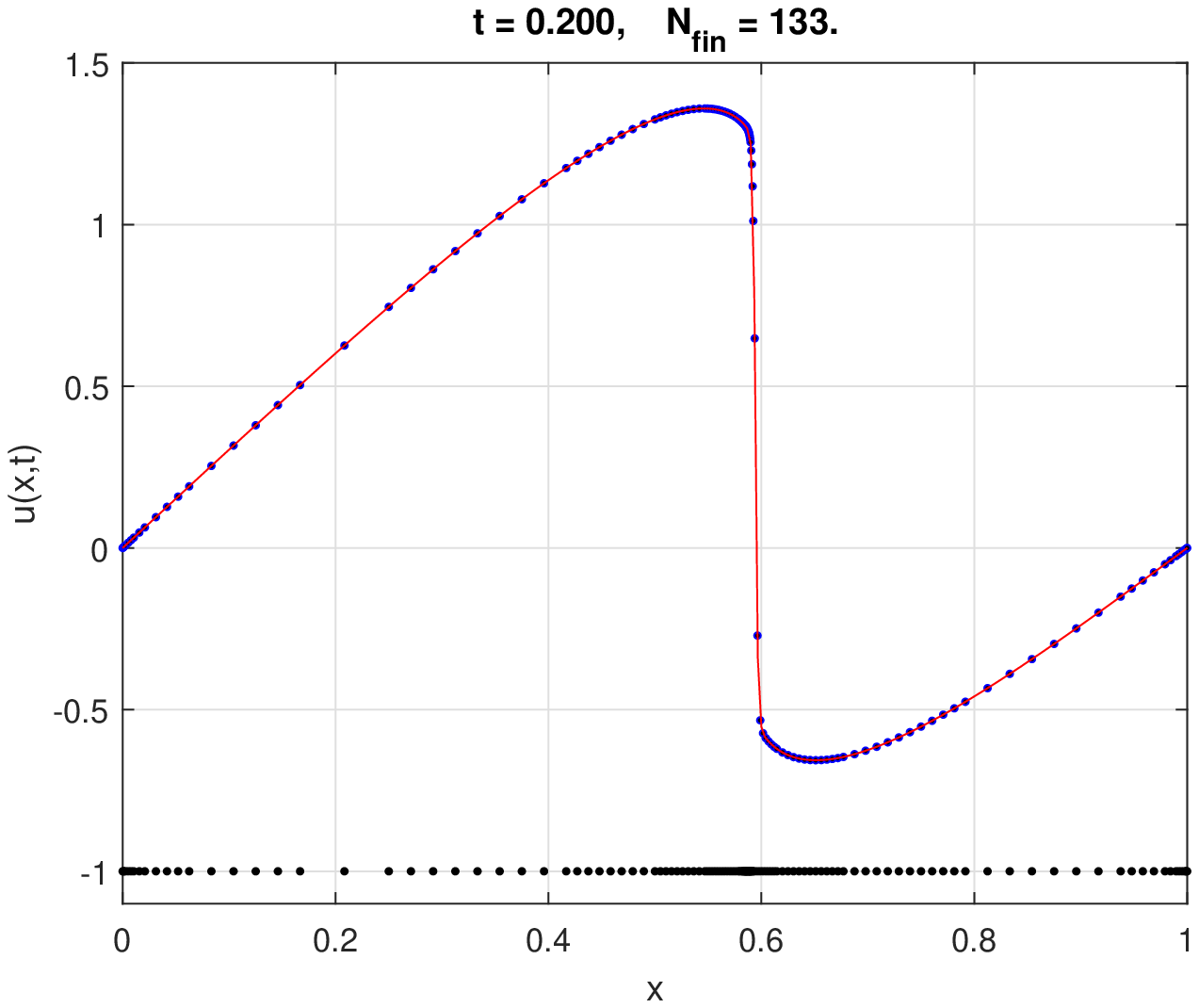}
\includegraphics[scale=0.26]{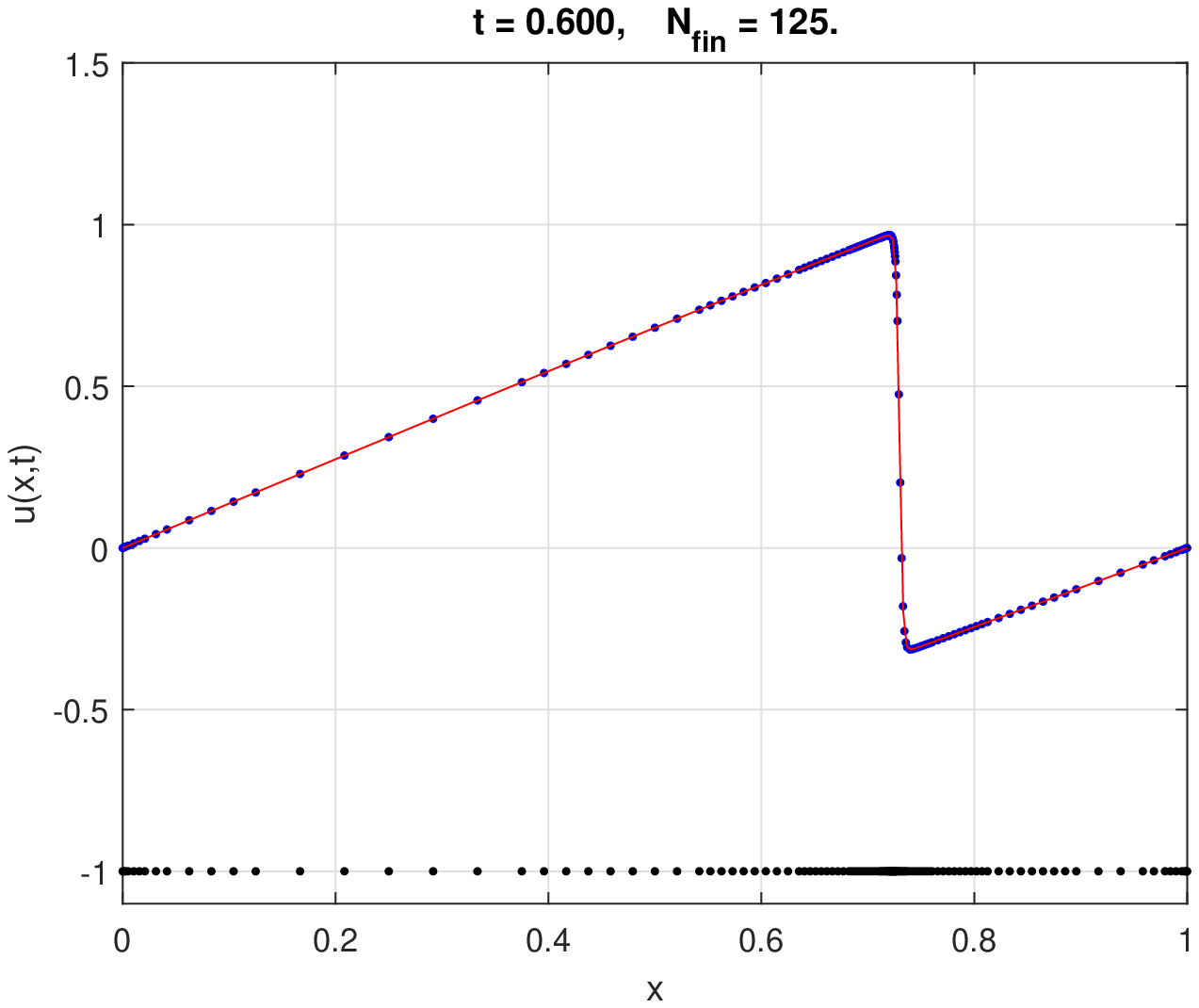} 
\includegraphics[scale=0.26]{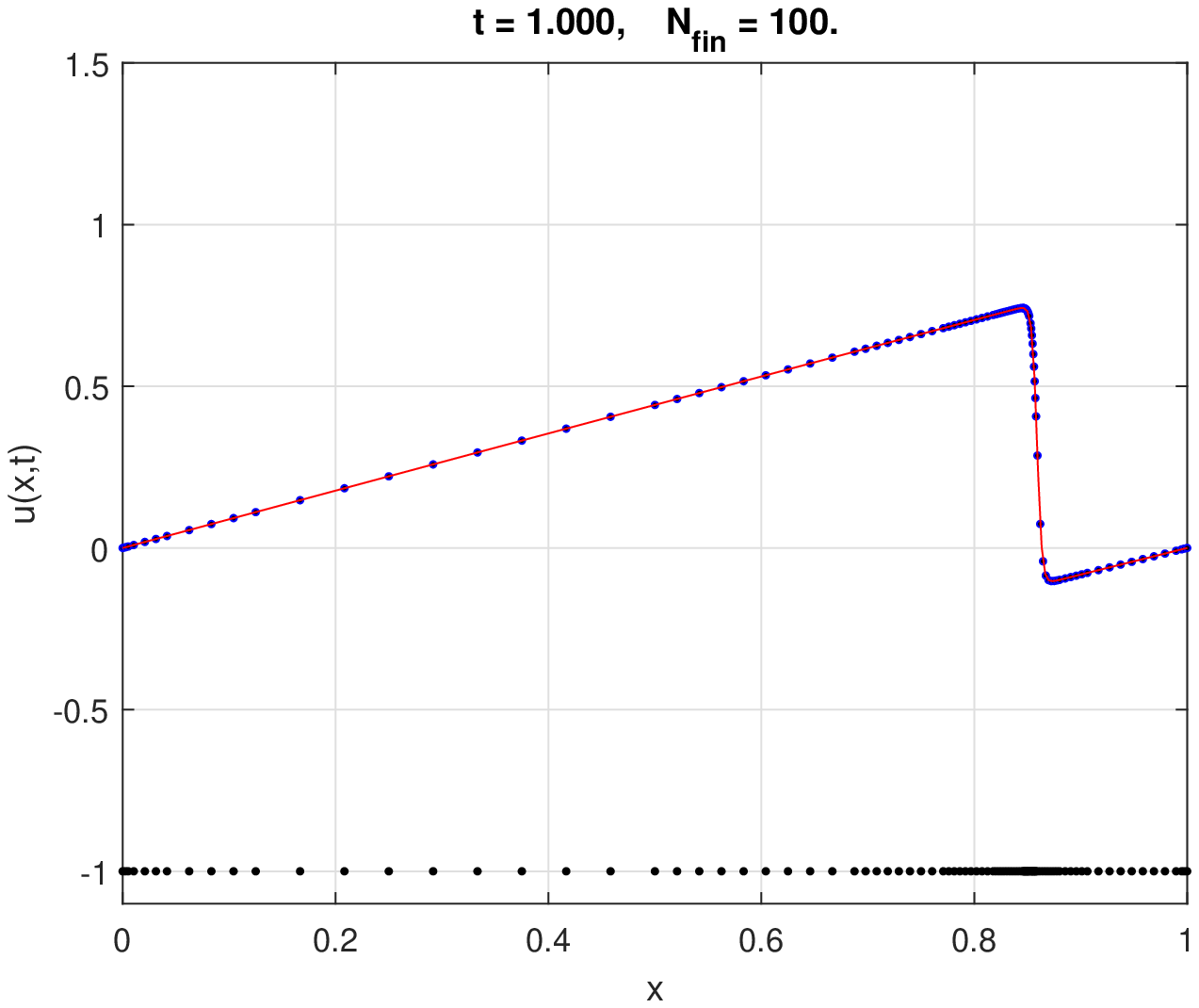}
}
\caption{Approximate solutions for Burgers equation of Example 2: $t=0$, $t=0.2$, $t=0.6$ and $t=1$ (left to right).}
\label{fig:3}
\end{figure}

Assuming in our experiments to fix the parameters $\nu = 10^{-3}$ and $\varepsilon_0=0.75$, we test our refinement scheme by starting initially from $N = 13$ equally-spaced collocation points. The adaptive LOOCV-based method is therefore applied for solving the Burgers problem at every time level until the final time arrives. In Figure \ref{fig:3} we show the solutions along with the final points obtained at time $t=0, 0.2, 0.6, 1$, with a threshold value $\tau=10^{-3}$. The related execution times are respectively given by CPU($t$) = $0.13, 1.54, 1.40, 1.10$. In general, we can notice that the adaptive performs well in each of the considered situations. These considerations can be extended even for other (more demanding) choices of the threshold $\tau$, but for the sake of shortness we omit such results.


\subsection{Case study 2: Solution of Allen-Cahn equations} \label{sec:42}

\begin{figure}[ht!]
\centering
{
\includegraphics[scale=0.26]{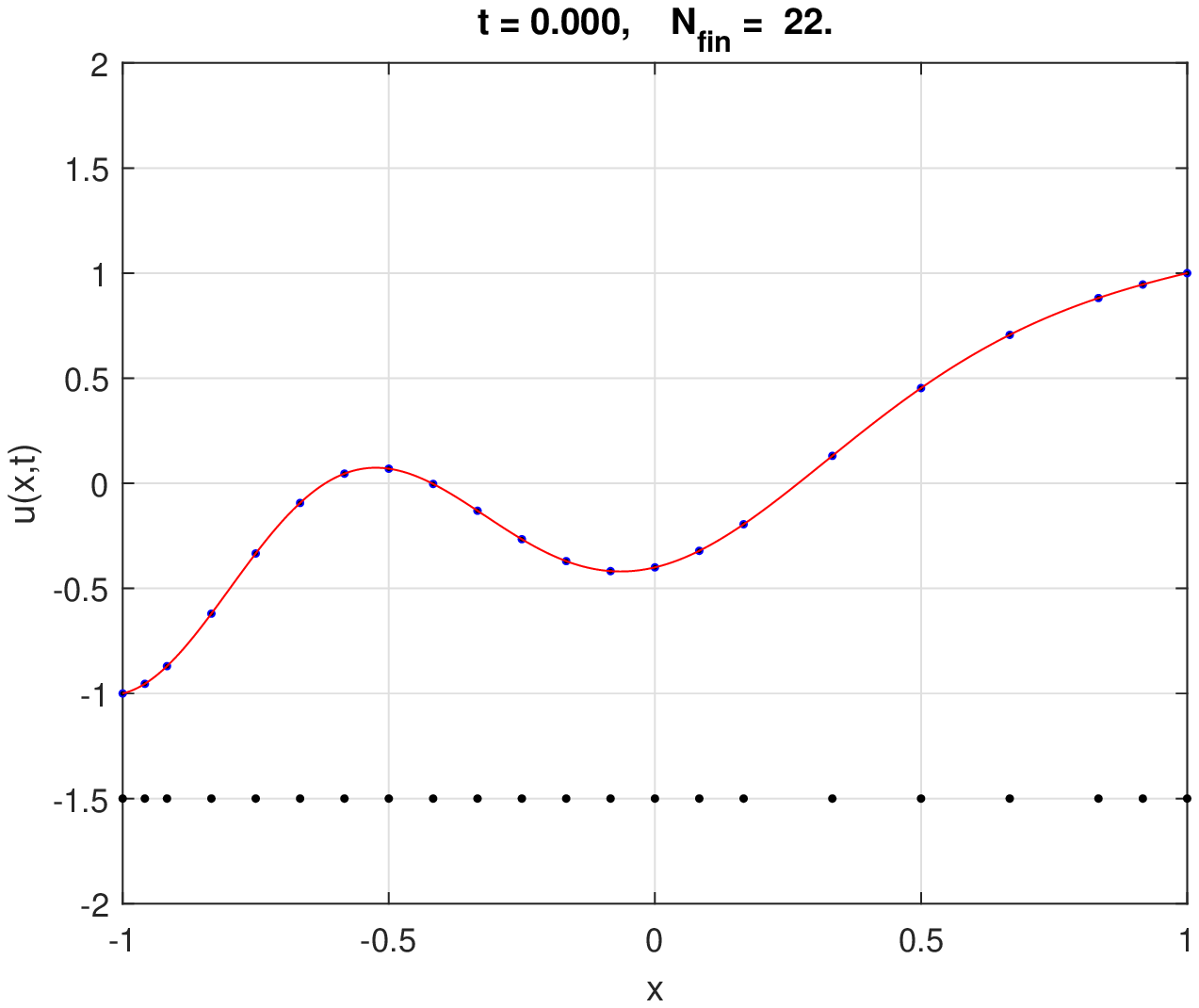} 
\includegraphics[scale=0.26]{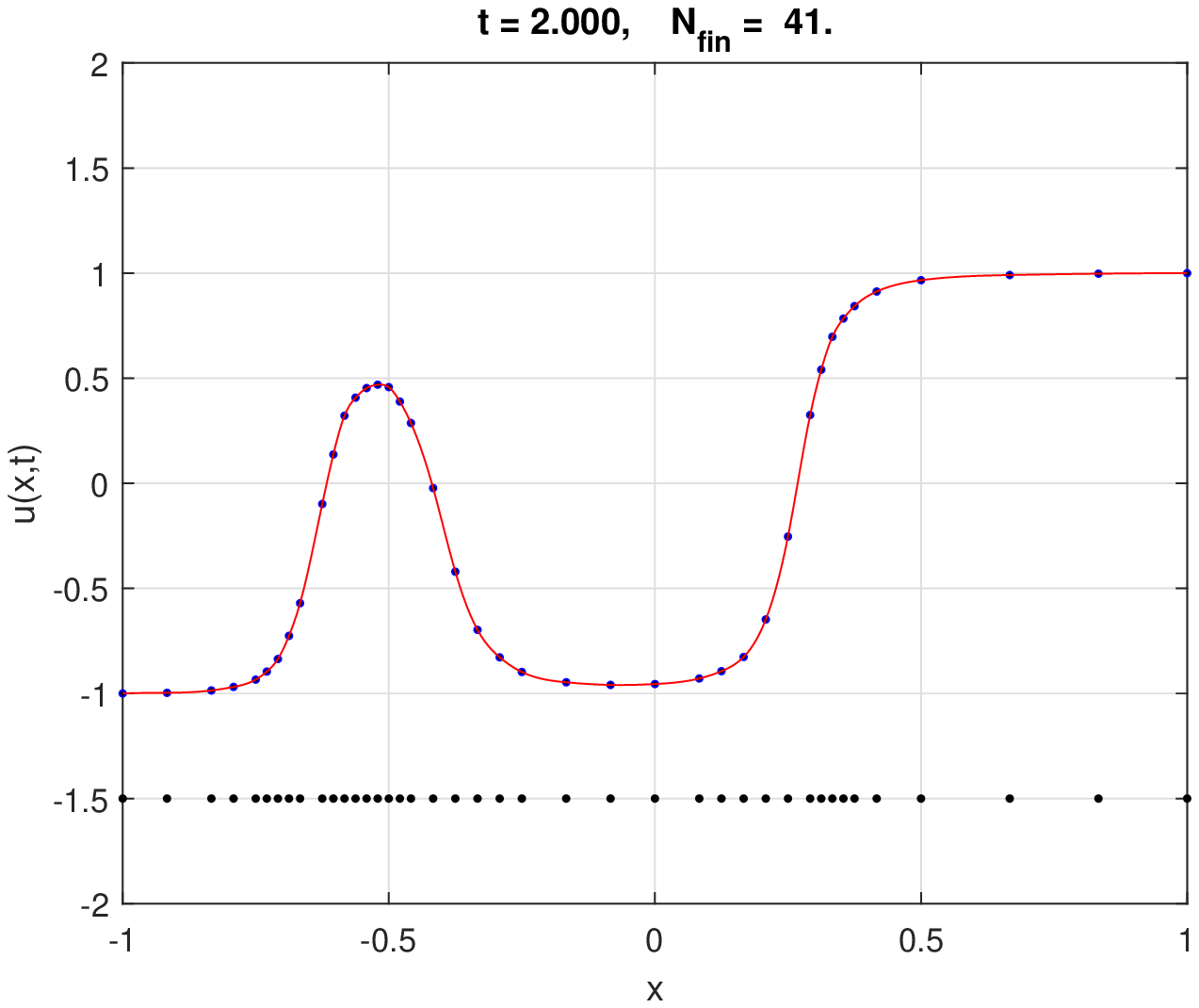} 
\includegraphics[scale=0.26]{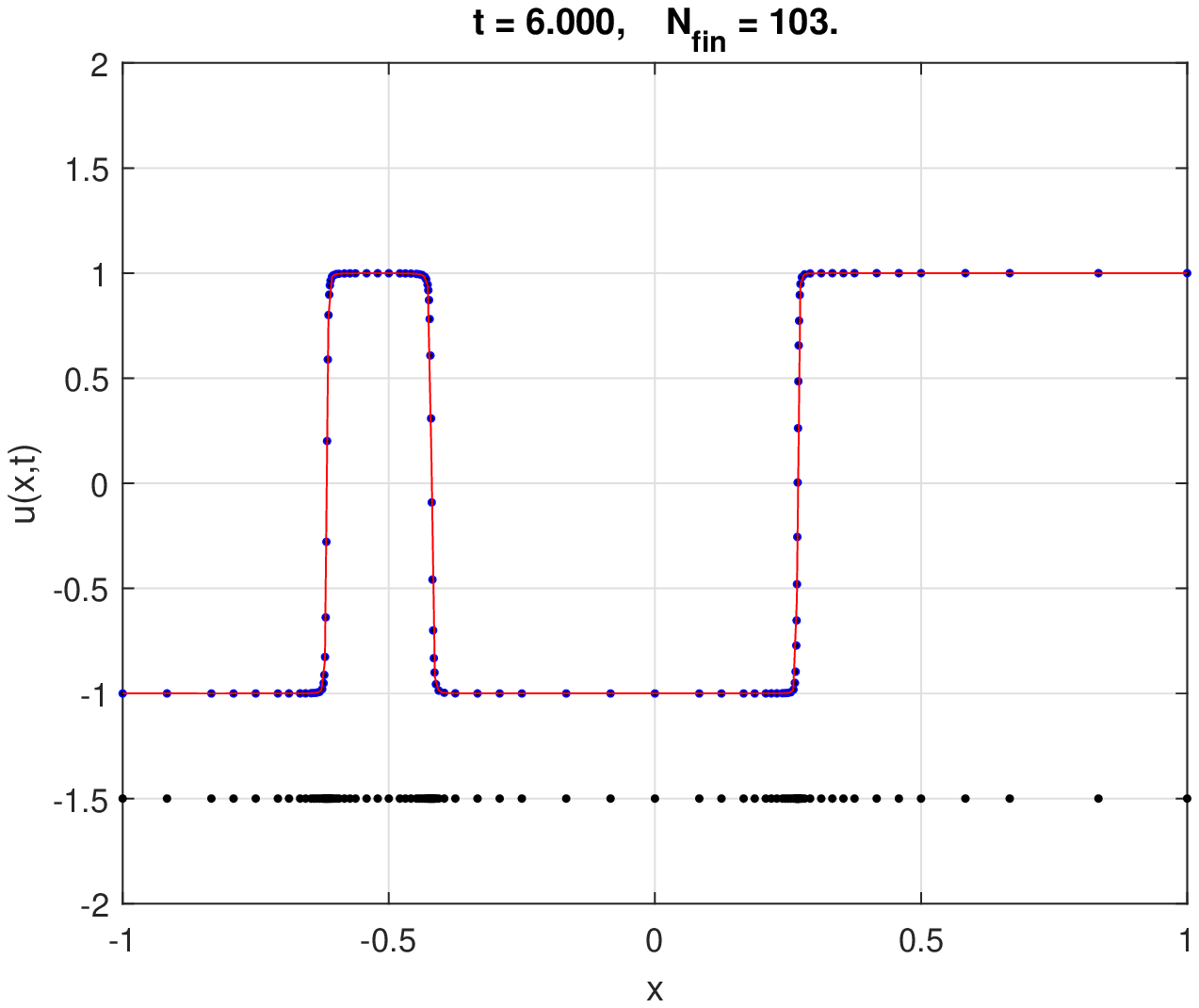} 
\includegraphics[scale=0.26]{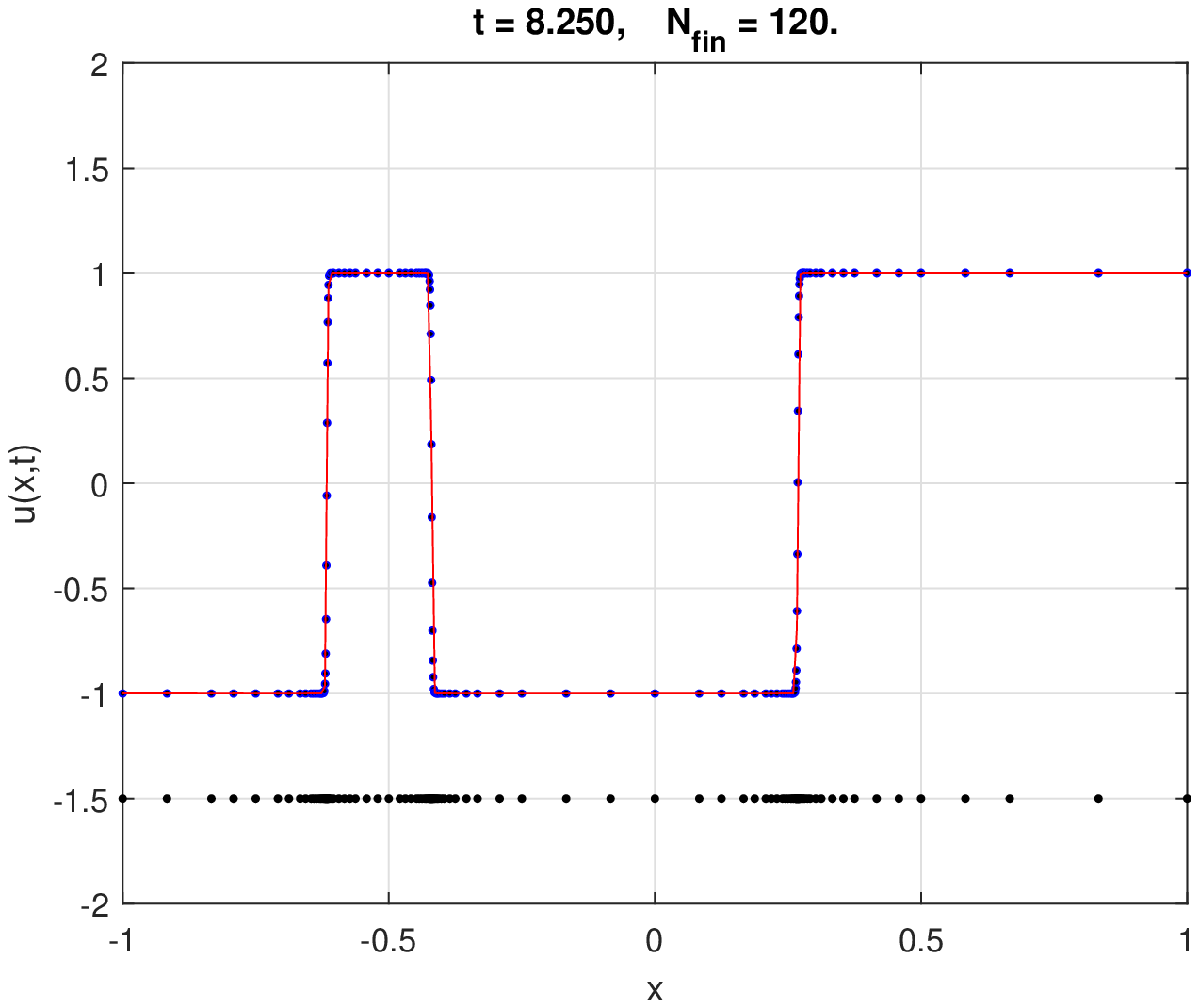}
}
\caption{Approximate solutions for Allen-Cahn equation: $t=0$, $t=2$, $t=6$ and $t=8.25$ (left to right). Such graphs refer to results given in Table \ref{tab:4}.}
\label{fig:4}
\end{figure}

In this subsection we consider another time-dependent BVP. In particular, referring to the generic problem \eqref{PDE-BC} defined on the computational domain $[a,b]\times [t_0,T] = [-1,1]\times [0,8.25]$, we now focus on the nonlinear Allen-Cahn equation \cite{dri07} given by
\begin{align} \label{allencahn}
\frac{\partial}{\partial t} u(x,t) = \nu \frac{\partial^2}{\partial x^2}u(x,t) + u(x,t)(1+u^2(x,t)), \qquad x \in (-1,1), \quad t \in (0,8.25],
\end{align}
where $\nu$ is a positive parameter. The initial condition is
\begin{align*} 
u(x,0) = 0.6x + 0.4 \sin\left[\frac{\pi}{2}\left(x^2-3x-1\right)\right],
\end{align*}
and the boundary conditions are 
\begin{align} \label{bc_allencahn}
u(-1,t) = -1, \quad u(1,t) = 1, \qquad t\in [0,8.25].
\end{align}
By setting the parameters $\nu = 10^{-6}$ and $\varepsilon_0=3$, the adaptation process is applied at every time step with threshold $\tau = 5\times 10^{-2}$, since the solution does not change its profile rapidly. Thus, the time interval $[0,8.25]$ is split by $M=34$ temporal levels. The algorithm starts initially with $N=13$ equally-spaced points in the spatial domain $[-1,1]$. In Table \ref{tab:4} we show the numerical results obtained for Allen-Cahn equation at time $t=0, 2, 4, 6, 8, 8.25$. From these experiments we can see that the number of collocation points required to achieve the stop criterion $\tau$ gradually grows with time $t$. In fact, we begin at time $t=0$ with an adaptive refinement that brings to $N_{fin}=22$ points until a maximum of $N_{fin}=120$ at $t=8.25$. However, the adaptive method turns out to be able to determine the steep variations of the solutions and is also sensitive to several oscillations. In Figure \ref{fig:4} we give some graphical representations of the found solutions together with the final sets of points for the above-mentioned values of $t$.

\begin{table}[ht!]
{\small
\begin{center} 
{
\begin{tabular}{cccccc} 
\midrule
 $t$ & $N_{fin}$  & CN($t$) & CPU($t$)   \\ \midrule
 $0.00$   & $22$  & $1.2$e${+3}$ & $0.07$  \\
 $2.00$   & $41$  & $6.7$e${+3}$ & $0.46$  \\
 $4.00$   & $76$  & $5.6$e${+4}$ & $0.46$  \\
 $6.00$   & $103$  & $3.2$e${+5}$ & $0.63$  \\
 $8.00$   & $120$  & $3.6$e${+5}$ & $0.43$  \\
 $8.25$   & $120$  & $3.6$e${+5}$ & $0.41$  \\
\midrule
\end{tabular}
}
\caption{Results obtained for Allen-Cahn equation by using the threshold value $\tau = 5\times 10^{-2}$ at some time levels $t$.} 
\label{tab:4}
\end{center}
}
\end{table}

\section*{Acknowledgments}
This work was partially supported by the INdAM-GNCS 2020 research project \lq\lq Multivariate approximation and functional equations for numerical modeling\rq\rq\ and by the 2020 project \lq\lq Mathematical methods in computational sciences\rq\rq\ funded by the Department of Mathematics \lq\lq Giuseppe Peano\rq\rq\ of the University of Torino. This research has been accomplished within the RITA \lq\lq Research ITalian network on Approximation\rq\rq\ and the UMI Group TAA \lq\lq Approximation Theory and Applications\rq\rq.






\begin{thebibliography}{99}







\bibitem{bur48} J. M. Burgers, A mathematical model illustrating the theory of turbulence, in: R. von Mises, T. von Karman (eds.), Advances in Applied Mechanics, Academic Press Inc., New York, 1948, pp. 171--199.












\bibitem{cav20a} R. Cavoretto, A. De Rossi, A two-stage adaptive scheme based on RBF collocation for solving elliptic PDEs, Comput. Math. Appl. 79 (2020) 3206--3222.

\bibitem{cav20b} R. Cavoretto, A. De Rossi, An adaptive LOOCV-based refinement scheme for RBF collocation methods over irregular domains, Appl. Math. Lett. 103 (2020) 106178.



\bibitem{cav21b} R. Cavoretto, Adaptive radial basis function partition of unity interpolation: A bivariate algorithm for unstructured data, J. Sci. Comput. 87 (2021) 41.

\bibitem{cav22} R. Cavoretto, A. De Rossi, A. Sommariva, M. Vianello, RBFCUB: A numerical package for near-optimal meshless cubature on general polygons, Appl. Math. Lett. 125 (2022) 107704.



\bibitem{che07} C.S. Chen, Y.C. Hon, R.A. Schaback, Scientific Computing with Radial Basis Functions, 2007.





\bibitem{dri07} T.A. Driscoll, A.R.H. Heryudono, Adaptive residual subsampling methods for radial basis function interpolation and collocation problems, Comput. Math. Appl. 53 (2007) 927--939.




\bibitem{fas07b} G.E. Fasshauer, J.G. Zhang, On choosing \lq\lq optimal\rq\rq\ shape parameters for RBF approximation, Numer. Algorithms 45 (2007) 345--368.


\bibitem{fas15} G.E. Fasshauer, M.J. McCourt, Kernel-based Approximation Methods using \textsc{Matlab}, Interdisciplinary Mathematical Sciences, Vol. 19, World Scientific Publishing Co., Singapore, 2015.









\bibitem{guo16} Y. Guo, Y. Shi, Y. Li, A fifth-order finite volume weighted compact scheme for solving one-dimensional Burgers' equation, Appl. Math. Comput. 281 (2016) 172--185.


\bibitem{har96} S.E. Harris, Sonic shocks governed by the modified Burgers' equation, European J. Appl. Math. 7 (1996) 201--222. 




\bibitem{hua94} W. Huang, Y. Ren, R.D. Russell, Moving mesh methods based on moving mesh partial differential equations, J. Comput. Phys. 113 (1994) 279--290.















\bibitem{lin21} L. Ling, F. Marchetti, A stochastic extended Rippa's algorithm for LpOCV,  Appl. Math. Lett. 129 (2022) 107955.


\bibitem{mar21} F. Marchetti, The extension of Rippa's algorithm beyond LOOCV, Appl. Math. Lett. 120 (2021) 107262.





\bibitem{naq17} S.L. Naqvi, J. Levesley, S. Ali, Adaptive radial basis function for time dependent partial differential equations, J. Prime Res. Math. 13 (2017) 90--106.




\bibitem{rip99} S. Rippa, An algorithm for selecting a good value for the parameter $c$ in radial basis function interpolation, Adv. Comput. Math. 11 (1999) 193--210.

















\bibitem{yan19} X. Yang, Y. Ge, L. Zhang, A class of high-order compact difference schemes for solving the Burgers' equations, Appl. Math. Comput. 358 (2019) 394--417.

\bibitem{yao15} G. Yao, J. Duo, C. S. Chen, L. H. Shen, Implicit local radial basis function interpolations based on function values, Appl. Math. Comput. 265 (2015) 91--102.



\end{thebibliography}
\end{document}